\documentclass[11pt]{amsart}
\usepackage{graphicx}
\usepackage{graphics}
\usepackage{amsmath}
\usepackage{amscd}
\usepackage{latexsym}

\begin{document}

\textwidth 6.2in
\textheight 7.6in
\evensidemargin .75in
\oddsidemargin.75in

\newtheorem{Thm}{Theorem}
\newtheorem{Lem}[Thm]{Lemma}
\newtheorem{Cor}[Thm]{Corollary}
\newtheorem{Prop}[Thm]{Proposition}
\newtheorem{Rm}{Remark}

\def\a{{\mathbb a}}
\def\C{{\mathbb C}}
\def\A{{\mathbb A}}
\def\B{{\mathbb B}}
\def\D{{\mathbb D}}
\def\E{{\mathbb E}}
\def\R{{\mathbb R}}
\def\P{{\mathbb P}}
\def\S{{\mathbb S}}
\def\Z{{\mathbb Z}}
\def\O{{\mathbb O}}
\def\H{{\mathbb H}}
\def\V{{\mathbb V}}
\def\Q{{\mathbb Q}}
\def\Cn{${\mathcal C}_n$}
\def\CM{\mathcal M}
\def\CG{\mathcal G}
\def\CH{\mathcal H}
\def\CT{\mathcal T}
\def\CF{\mathcal F}
\def\CA{\mathcal A}
\def\CB{\mathcal B}
\def\CD{\mathcal D}
\def\CP{\mathcal P}
\def\CS{\mathcal S}
\def\CZ{\mathcal Z}
\def\CE{\mathcal E}
\def\CL{\mathcal L}
\def\CV{\mathcal V}
\def\CW{\mathcal W}
\def\IC{\mathbb C}
\def\IF{\mathbb F}
\def\IK{\mathcal K}
\def\IL{\mathcal L}
\def\IP{\bf P}
\def\IR{\mathbb R}
\def\IZ{\mathbb Z}

\title{The Akhmedov-Park exotic  ${\C\P}^{2}\#2\bar{\C\P}^2$}
\author{Selman Akbulut}
\thanks{The author is partially supported by NSF grant DMS 0905917}
\keywords{}
\address{Department  of Mathematics, Michigan State University,  MI, 48824}
\email{akbulut@math.msu.edu }
\subjclass{58D27,  58A05, 57R65}
\date{\today}
\begin{abstract} 
Here we draw a handlebody picture for the exotic ${\C\P}^{2}\#2\bar{\C\P}^2$ constructed by Akhmedov and Park.

\end{abstract}

\date{}
\maketitle

\setcounter{section}{-1}

\vspace{-.34in}

\section{Introduction}

Akhmedov-Park manifold  $N$ \cite{ap1} is an irreducible symplectic $4$-manifold which is an exotic copy of ${\C\P}^{2}\#2\bar{\C\P}^2$. It is a variant of the exotic ${\C\P}^{2}\#3\bar{\C\P}^2$ which they constructed previously \cite{ap2} (with its  handlebody  in \cite{a1}). $N$ consists of  two codimension zero pieces  glued along the common boundaries: $$N=\tilde{E}_{0}\smile_{\partial} \tilde{E}_{1}$$  
The two pieces are constructed as follows:  We start with the product of a genus $2$ surface and the torus $E = \Sigma_{2}\times T^{2}$. Let $<a_1,b_1,a_2,b_2>$ and $<C,D >$ be the standard circles generating the first homology of $\Sigma_{2}$ and $T^{2}$ respectively (the cores of the $1$-handles). Then $\tilde{E}$ is obtained from $E$ by doing ``Luttinger surgeries'' to the four subtori $(a_1\times C, a_1)$, $(b_1\times C, b_1)$,   $(a_2 \times C, C)$  $(a_2 \times D, D)$ (see \cite{a1} for a handle description of $\tilde{E}$).  Then
$$\tilde{E_{0}}=\tilde{E}- \Sigma_{2}\times D^2$$

To construct the other piece, let $K\subset S^3$ be the trefoil knot, and $S^{3}_{0}(K)$ be the $3$-manifold obtained doing $0$-surgery to $S^3$ along $K$. Being a fibered knot, $K$ induces a fibration $T^2\to S^{3}_{0}(K)\to S^1$ and the fibration
$$T^2\to S^{3}_{0}(K)\times S^1\to T^2$$
Let $T^{2}_1$ and $T^{2}_2$ be the vertical (fiber) and the horizontal (section) tori of this fibration, transversally intersecting one point $p$. The degree $2$ map $S^1\to S^1$ induces an imbedding $S^1\times D^3 \hookrightarrow   S^1\times D^3$, and by using one of the circle factors of $T^2$  it induces an imbedding $\rho: T^2\times D^2 \hookrightarrow T_{1}^2\times D^2$. The torus $T^{2}_{3}:= \rho(T^2\times 0)\subset T^{2}_{1}\times D^2$  intersects a normal disk  $D^2$ of $T_{1}^2$ transversally at two points $\{p_1,p_2\}$, hence it intersects $T^{2}_{2}$ transversally at  $\{p_1,p_2\}$. Therefore  by blowing up the total space at $p_1$, we get a pair of imbedded pair tori 
$T^{2}_{3}\cup T^{2}_{2} \subset (S^{3}_{0}(K)\times S^1) \# \bar{\C\P}^2$, each with self intersection $-1$, and transversally intersecting each other at one point $p_2$. Smoothing $T_{2}^{2}\cup T^{3}_{2}$ at $p_2$ gives a genus $2$ surface $\Sigma_{2} \subset (S^{3}_{0}(K)\times S^1) \# \bar{\C\P}^2$ with trivial normal bundle $\Sigma_{2}\times D^2$, then we define: 

$$\tilde{E}_{1}=(S^{3}_{0}(K)\times S^1) \#  \bar{\C\P}^2 - \Sigma_{2} \times D^2$$


\noindent In short  $M= \tilde{E}\;\sharp_{\;\Sigma_{2}} (S^{3}_{0}(K)\times S^1) \#  \bar{\C\P}^2 $ ($\sharp_{\;\Sigma_{2}}$ denotes fiber sum along $\Sigma_{2}$). There is a related construction of  exotic ${\C\P}^{2}\#2\bar{\C\P}^2$ given in \cite{fs}, which appeared after \cite{ap1} (in fact  arxiv.org version of \cite{ap1} appeared in 2007).  We thank Anar Akhmedov for describing this beautiful manifold to us.

\section{construction of $N$}

In \cite{a1} a handlebody picture for $\tilde{E}_{0}$ is already given. Recall, the starting point was  Figure 1  from  \cite{a2}, which is a handelbody of 
$E_{0}= E- \Sigma_{2}\times D^2$ and the desrcription of a diffeomorphism $f : \partial E_{0}\to \partial (\Sigma_{2}\times D^2)$. In \cite{a1} by using this  we constructed a handlebody for $\tilde{E}_{0}$ and gave a description of a diffeomorphism $f : \partial \tilde{E}_{0}\to \partial (\Sigma_{2}\times D^2)$ (Figure 2). In Figure 2 we also indicate the generators of $\pi_{1}(\tilde{E}_{0})$. Next we will give a handlebody for $\tilde{E}_1$.

\vspace{.1in}

The first picture of Figure 3 is the handlebody of $S^{3}_{0}(K) \times S^{1}$ (\cite{a3}), where in this picture the horizontal $T_{2}\times D^2$ is clearly visible, but not the vertical torus $T^{2}_{1}$ (which consists of the Seifert surface of $K$ capped off by the $2$-handle bounded by the trefoil knot). In the second picture of Figure 3 we redraw this handlebody so that both the vertical and the horizontal tori are clearly visible (reader can check this by canceling $1$- and $2$- handle pairs from the second picture to obtain the first picture). By an isotopy and creating a canceling $1$- and $2$- handle pair we obtain the first and the second pictures of Figure 4.  Now we want to find an imbedded copy of $T^2\times D^2$ inside of the neighborhood  $T_{1}^{2}\times D^2$ of the vertical torus (indicated by dotted rectangle) imbedded by the degree two map, i.e. we want to locate  $T^{2}_{3}$. 

\vspace{.1in}

Reader can bemuse herself by first trying to draw a handle description of a copy of $S^1\times D^3$ (a circle with dot) inside $S^1\times D^3$ (another circle with dot) imbedded by a degree two map. After several tries it will become clear that one should first start with $S^1\times D^2$ inside of another $S^1\times D^2$ imbedded by a degree two map. This is given by the first picture of Figure 5 (Heegard picture). Here the $1$-handle of the sub $S^1\times D^2$ is $B$ and the $1$-handle of the ambient $S^1\times D^2$ is $A$. Crossing this picture $S^1$ we get a handle description of $T^2\times D^2$  inside of another $T^{2}\times D^2$ imbedded by the degree two map, which is the second picture of Figure 5.  

\vspace{.1in}

By changing the $1$-handle presentation from pair of balls to circle with dot notation of \cite{a4}, we get the pictures of Figure 6, where the degree two imbedding of $T^{2}\times D^2$  inside of $T^{2}\times D^2$ is  now clearly  visible. For emphasis, we should mention that the either picture of Figure 6 is just $T^2\times D^2$, drawn in a nonstandard way so that the degree two imbedded copy of the sub $T^2\times D^2$  is clearly visible inside. Now we want to carry this picture to Figure 4 install it inside of the dotted rectange (i.e. inside of $T_{1}^{2}\times D^2$). 

\vspace{.1in} 

For this purpose, in Figures 6-8, we describe a diffeomorphism h (in terms of handle slides) from the copy of $T^2\times D^2$ in Figure 6 to the standard copy of $T^2\times D^2$, and then apply the inverse $h^{-1}$ of this diffeomorphism to the region enclosed by the dotted rectangle in Figure 4. This gives the second picture of Figure 9.  Reader can quickly  verify this by simply applying $h$ to the second picture of Figure 9 and get the first picture of Figure 9. Think of Figure 9 as a piece of $S^{3}_{0}(K)\times S^1$ (i.e. the region enclosed by dotted rectangle in Figure 4). Now we want to work  only in this region (i.e. manipulate it by diffeomorphisms), then at the end install it back in Figure 4. This will help to simplify our pictures, otherwise we have to carry the rest of whole Figure 4 every step, even though we are only working inside of $T_{1}^2\times D^2$.

\vspace{.1in} 

By blowing up  $S^{3}_{0}(K)\times S^1$ we go from Figure 9 to Figure 10,  which is $(S^{3}_{0}(K)\times S^1)\# \bar{\C\P}^2$ (actually the figure shows only piece of it). The thick arrow in Figure 9 indicates where we performed the blowing up operation. Then by isotopies, handle slides (indicated by dotted arrows), and handle cancelations we go from Figure 10 to Figure 14.  By installing Figure 14 inside of the dotted square in Figure 4 we obtain the whole picture of $(S^{3}_{0}(K)\times S^1)\# \bar{\C\P}^2$, which is Figure 15. Figure 16 is an equivalent way of expressing  Figure 15, where we draw the slice $1$-handle as a pair standard $1$-handles and a $2$-handle. Now finally, in Figures 15 and 16 we see clearly the imbedded copy  $\Sigma_{2}\times D^2 \subset (S^{3}_{0}(K)\times S^1) \# \bar{\C\P}^2$ discussed in the introduction. 

\vspace{.1in} 

To obtain $N$ we need to  remove this copy of  $\Sigma_{2} \times D^2$  from inside of $(S^{3}_{0}(K) \times S^{1} )\# \bar{\C\P}^2$ and replace it with $\tilde{E_{0}}$. The arcs in Figure 2 (describing the diffeomorphism $f$) show us how to do this, resulting with the handlebody picture of $M$ in Figure 17. Throughout the paper we use the convention that when we don't  indicate the framing a $2$-handle it means it is zero framed.

\section{ Triviality of $ \pi_{1}(N)$}

Recall that in \cite{a1}  the fundamental group $\pi_{1}(\tilde{E}_{0})$ is given in terms of the generators of Figure 2 (by reading off the relations by tracing the attaching knots of the $2$-handles, starting at the points indicated by small circles)

 \begin{equation*}
 \pi_{1}(\tilde{E_0}) = \left\langle a,\; c,\; e,\; g,\; h,\; q \; \left| 
\begin{array}{cc}
[a,e]=1,  \;  [h,e]=1,  \;  [a,q]=c, & \\
 $$[c^{-1},a]=cq, \; [g,a]=h, \; [g^{-1},h]=a,
  & \\
$$[cq,ah]=1,\; [c,ah]=1,\;  [q^{-1},c][g^{-1},e]=1&
\end{array}  \right\rangle \right.
\end{equation*}

Now we can calculate  $ \pi_{1}(N)$ from this presentation  by introducing the new generators $x,y,u,v,w$ of Figure 17,  and  introducing $11$ new relations coming from the $11$ new $2$-handles of Figure 17 (attached on top of $\tilde{E_0}$). These new relations can be easily red off from the Figure 17, they are given by the following words: 

$$vg^{-1}h^{-1}gha,\;  xyxy^{-1}x^{-1}y^{-1},\; cyx^{-1}c^{-1}xy^{-1},\; qy^{-1},\;  xyx^{-1}wy^{-1},$$
$$yx^{-1}v^{-1}u^{-1},\;vu^{-1}vu,\;uwu^{-1}vw^{-1},\;u^{2}e^{-1},\; we^{-1}w^{-1}e,\; wg^{-1}$$

After eliminating $y, w$ from the short words, these new relations become:

$$ v =a^{-1}[h^{-1},g^{-1}],\; xqx=qxq,\; cqx^{-1}c^{-1}xq^{-1}=1$$
$$ g=[x,q^{-1}],\; qx^{-1}=uv=v^{-1}u,\;,v=[u,g^{-1}],\; u^2=e,\; [g,e]=1 $$

\vspace{.1in}

The last relation of $\pi_{1}(\tilde{E}_{0})$  and $ [g,e]=1$  implies $ [q,c]=1$, and this  together with the third relation above implies $[c,x]=1$.  One can make a few more simplifications, for example in the presence of  $[q,c]=1$ the third relation above is equivalent to $[c,x]=1$. Also the $7$th and $8$th relations of $\pi_{1}(\tilde{E}_{0})$  together are equivalent to  $[q,ah]=1$ and $[c,ah]=1$. So we have:

\vspace{.05in}

   \begin{equation*}
 \pi_{1}(N) = \left\langle a,\; c,\;  g,\; h,\; q, \; x,\; u,\;v \left| 
\begin{array}{cc}
[a,u^2]=1,  \;  [h,u^2]=1,  & \\
 $$ [a,q]=c,\; [a^{-1},c^{-1}]=q, \; [g,a]=h
  & \\
$$[g^{-1},h]=a,\; [q,ah]=1,\; [c,ah]=1,\;  [q,c]=1&\\
v=a^{-1}[h^{-1},g^{-1}], \; xqx=qxq,\; [c,x]=1&\\
g=[x,q^{-1}],\; qx^{-1}=uv=v^{-1}u,&\\
  v=[u,g^{-1}],\; [g,u^2]=1
\end{array}   \hspace{-.15in} \right\rangle \right.
\end{equation*}

\vspace{.1in}

\noindent From this presentation by using the group theory software GAP, reader can verify that $\pi_{1}(N)$ is the trivial group. Also notice  $g=[x,q^{-1}]$ and $xqx=qxq$ $\Rightarrow$ $g^{-1}x^{2}=xq$; again this with $xqx=qxq$ gives an equivalent presentation $x=gqg^{-1}$. This helps to reduce number of generators. Also by GAP, we can check that nine of the above relations are redundant, and the group has the following presentation, and it is trivial, though it would be nice to find a by hand proof of the triviality of this group as in the case of \cite{a1}. 

 \begin{equation*}
 \pi_{1}(N) = \left\langle a,\; c,\;  g,\; h,\; q  \left| 
\begin{array}{cc}
 $$ [a,q]=c,\; [a^{-1},c^{-1}]=q, \; [g,a]=h,\; [g^{-1},h]=a
  & \\
$$[q,ah]=1,\; [c,ah]=1,\;  gq^{-2}g=q^{-1}gq^{-1}
\end{array}  \hspace{-.15in} \right\rangle \right. = \langle 1\rangle
\end{equation*}

We leave this and the topological  consequences of this handlebody as a subject of  a future paper. It is interesting to observe how the manifold $M$ of  \cite{a1} evolves into this manifold $N$. The pictures suggest that $N$ might be obtained by some sort of `rational blowing down" (and maybe mixed in with some knot surgery) operation to $M$.

    \begin{figure}[ht]  \begin{center}  
\includegraphics[width=.9\textwidth]{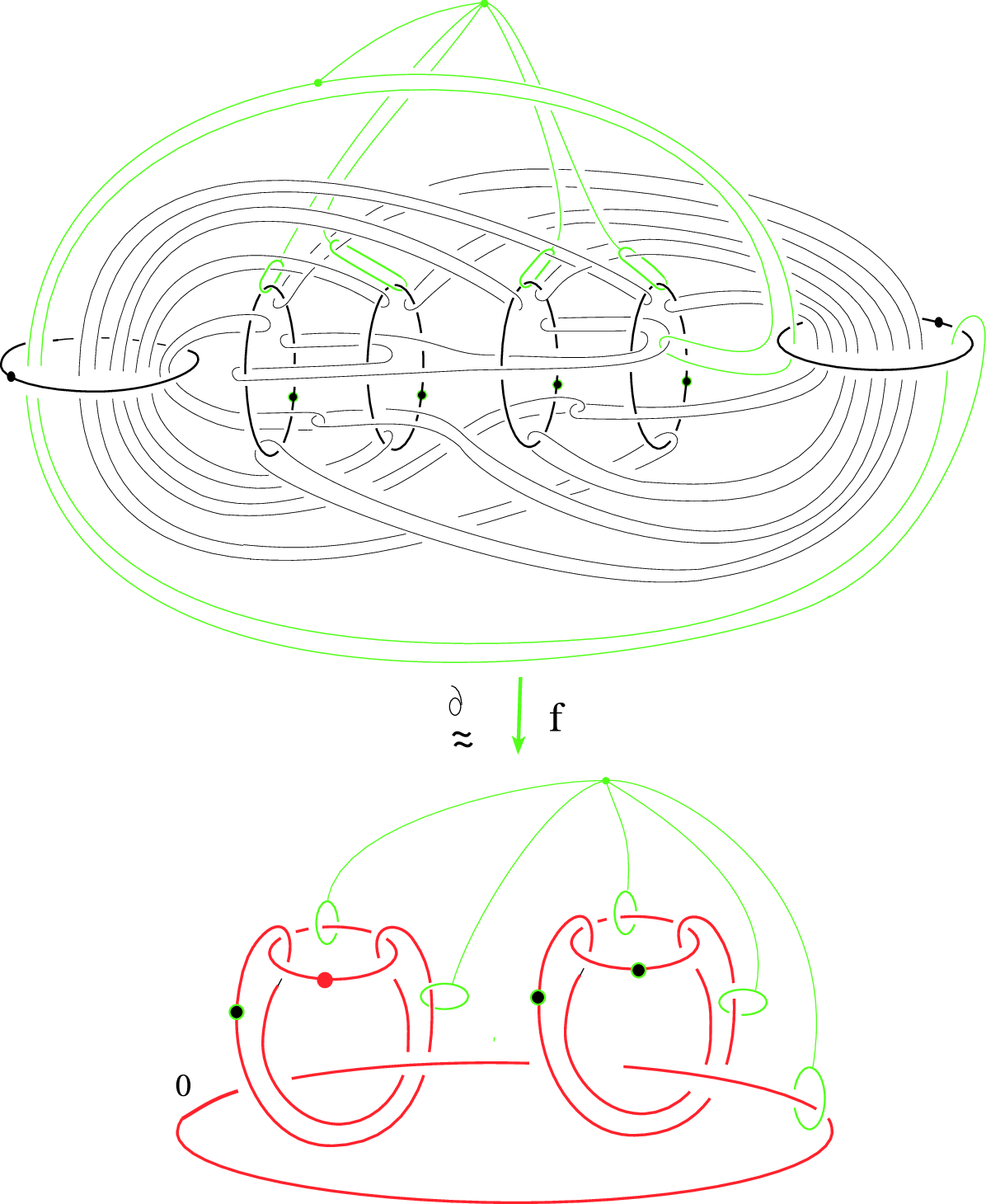}   
\caption{$E_{0}$ and  $f: \partial E_{0}\stackrel{\approx}{\longrightarrow} \Sigma_{2}\times S^1$} 
\end{center}
\end{figure} 

    \begin{figure}[ht]  \begin{center}  
\includegraphics[width=.95\textwidth]{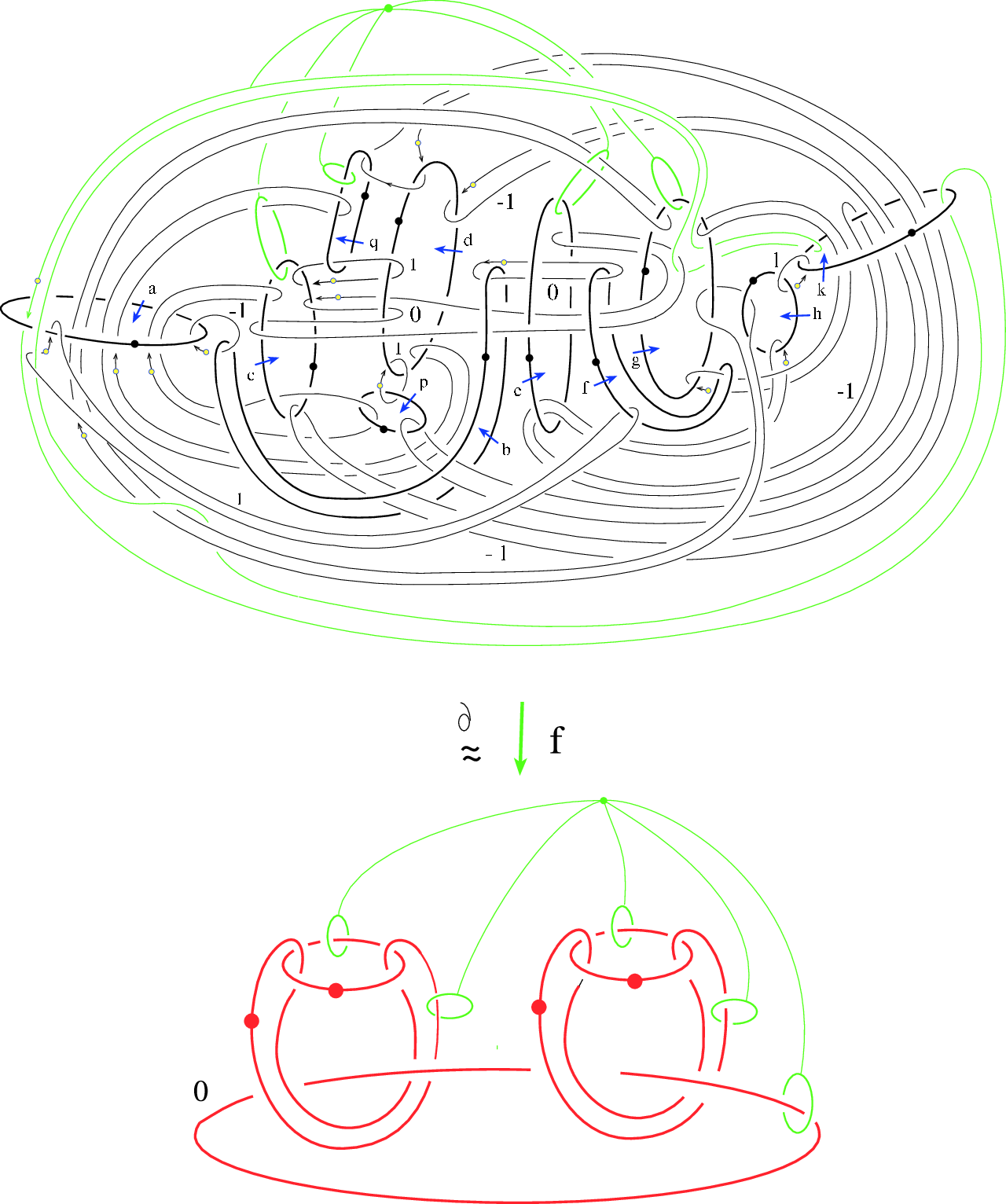}   
\caption{$\tilde{E}_{0}$ and  $f: \partial \tilde{E}_{0}\stackrel{\approx}{\longrightarrow} \Sigma_{2}\times S^1$ } 
\end{center}
\end{figure} 

    \begin{figure}[ht]  \begin{center}  
\includegraphics[width=.7\textwidth]{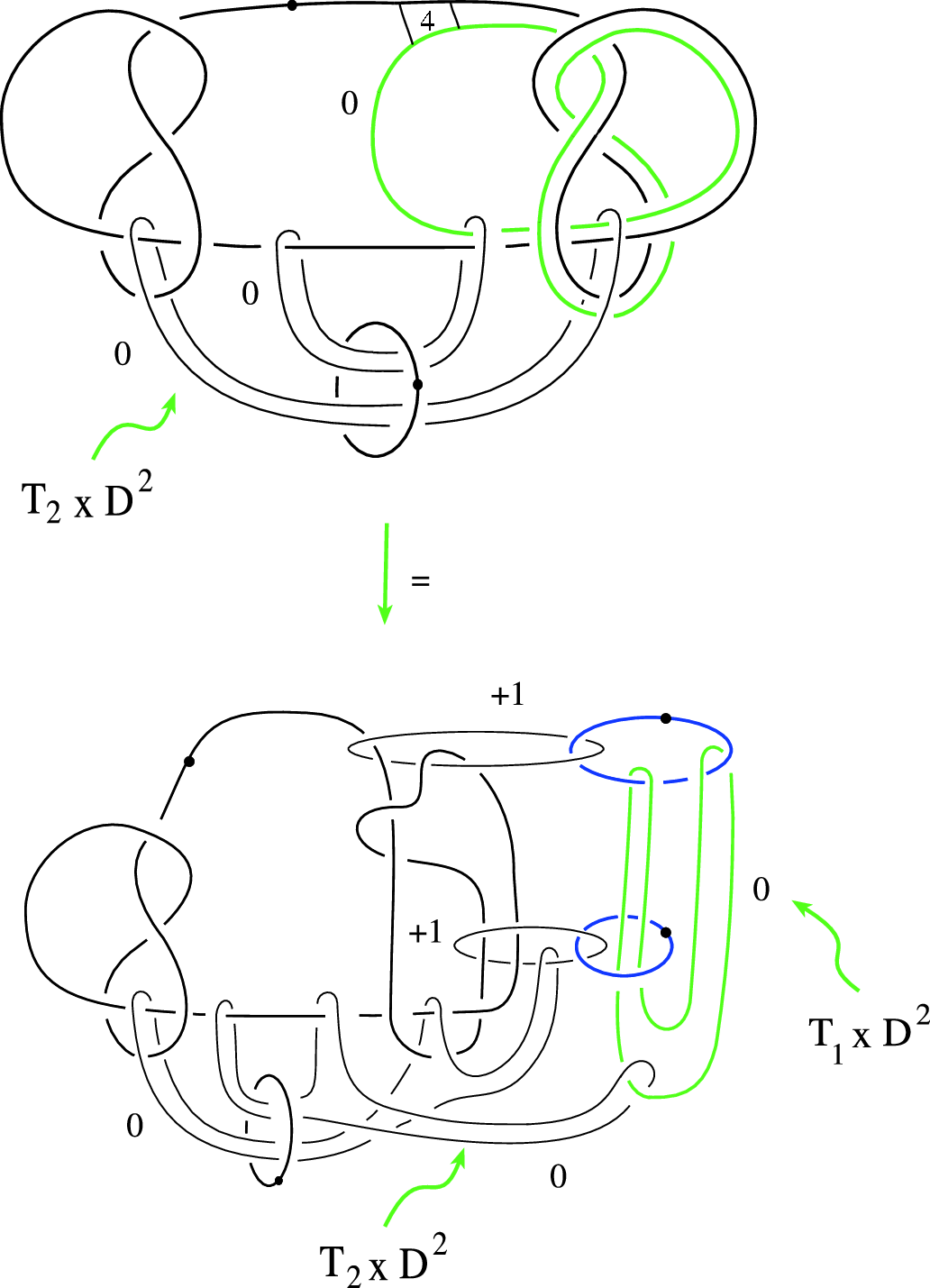}   
\caption{ $S^{3}_{0}(K)\times S^1$} 
\end{center}
\end{figure} 

    \begin{figure}[ht]  \begin{center}  
\includegraphics[width=.7\textwidth]{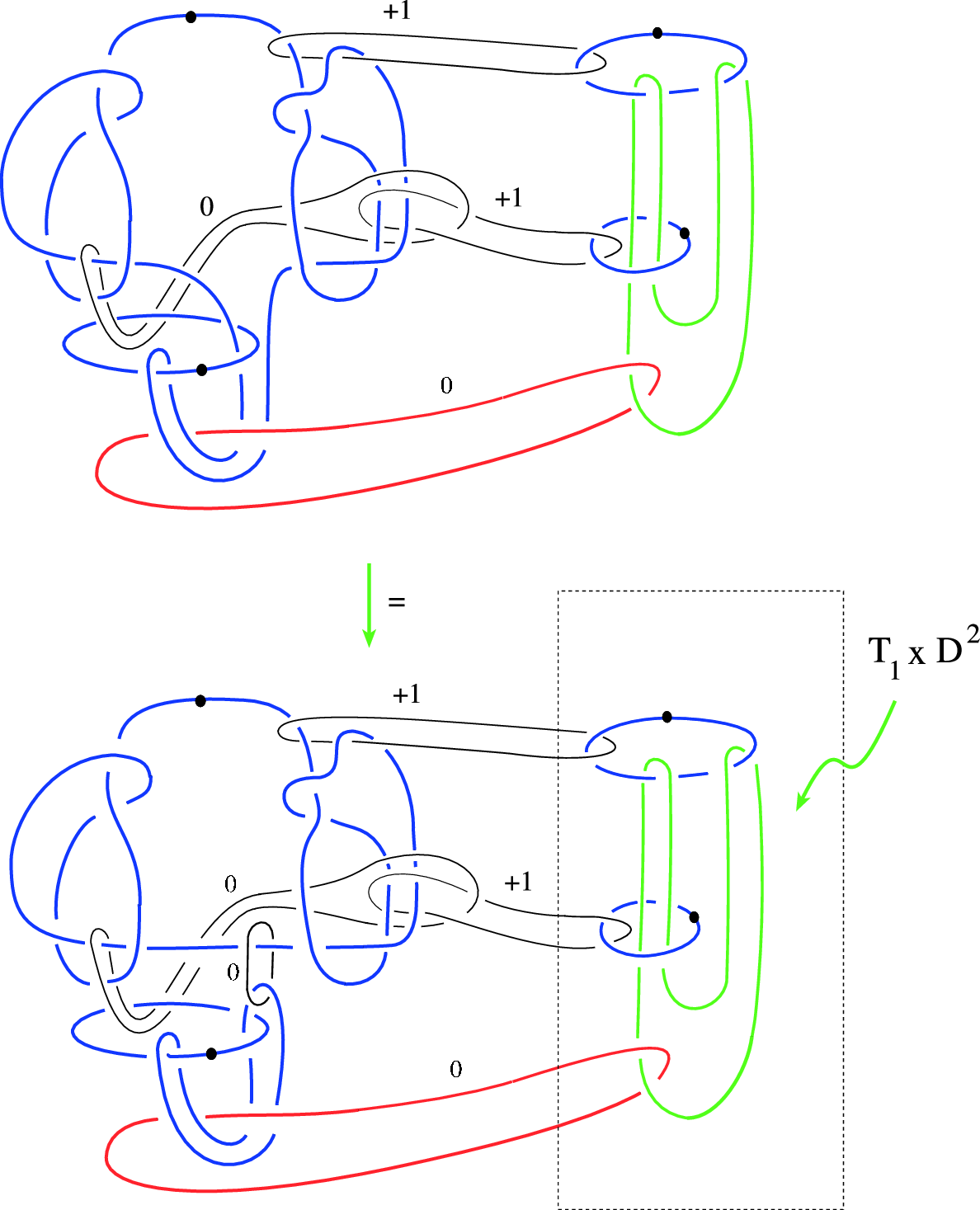}   
\caption{  $S^{3}_{0}(K)\times S^1$} 
\end{center}
\end{figure}

    \begin{figure}[ht]  \begin{center}  
\includegraphics[width=.6\textwidth]{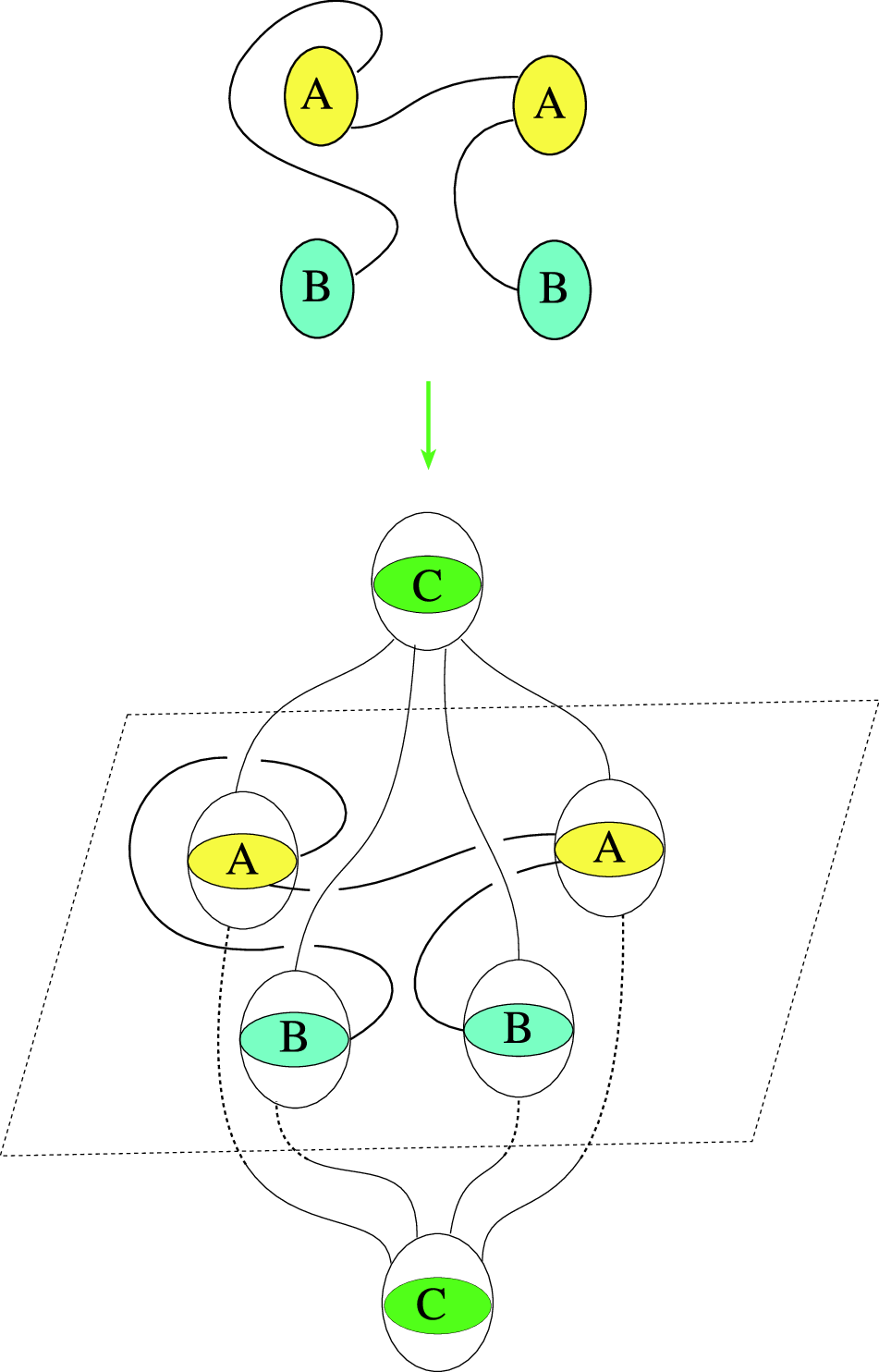}   
\caption{Building degree $2$ imbedding $T^{2}\times D^2 \hookrightarrow T^2\times D^2$ by crossing degree $2$ imbedding $S^{1}\times D^2 \hookrightarrow S^{1}\times D^2$ with $S^1$} 
\end{center}
\end{figure} 

    \begin{figure}[ht]  \begin{center}  
\includegraphics[width=.6\textwidth]{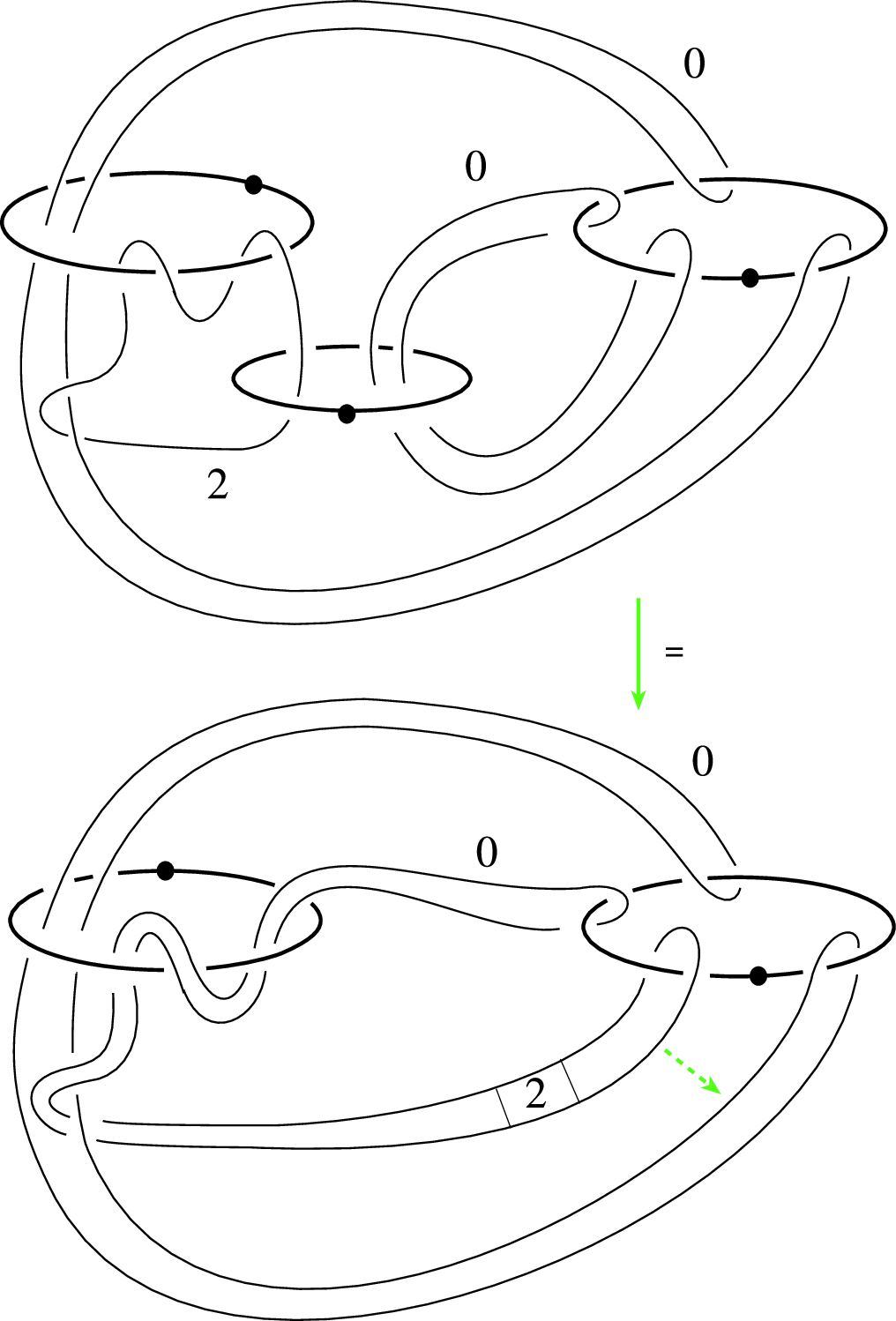}   
\caption{$T^2 \times D^2$ containing a degree $2$ imbedded copy of itself} 
\end{center}
\end{figure}

    \begin{figure}[ht]  \begin{center}  
\includegraphics[width=.7\textwidth]{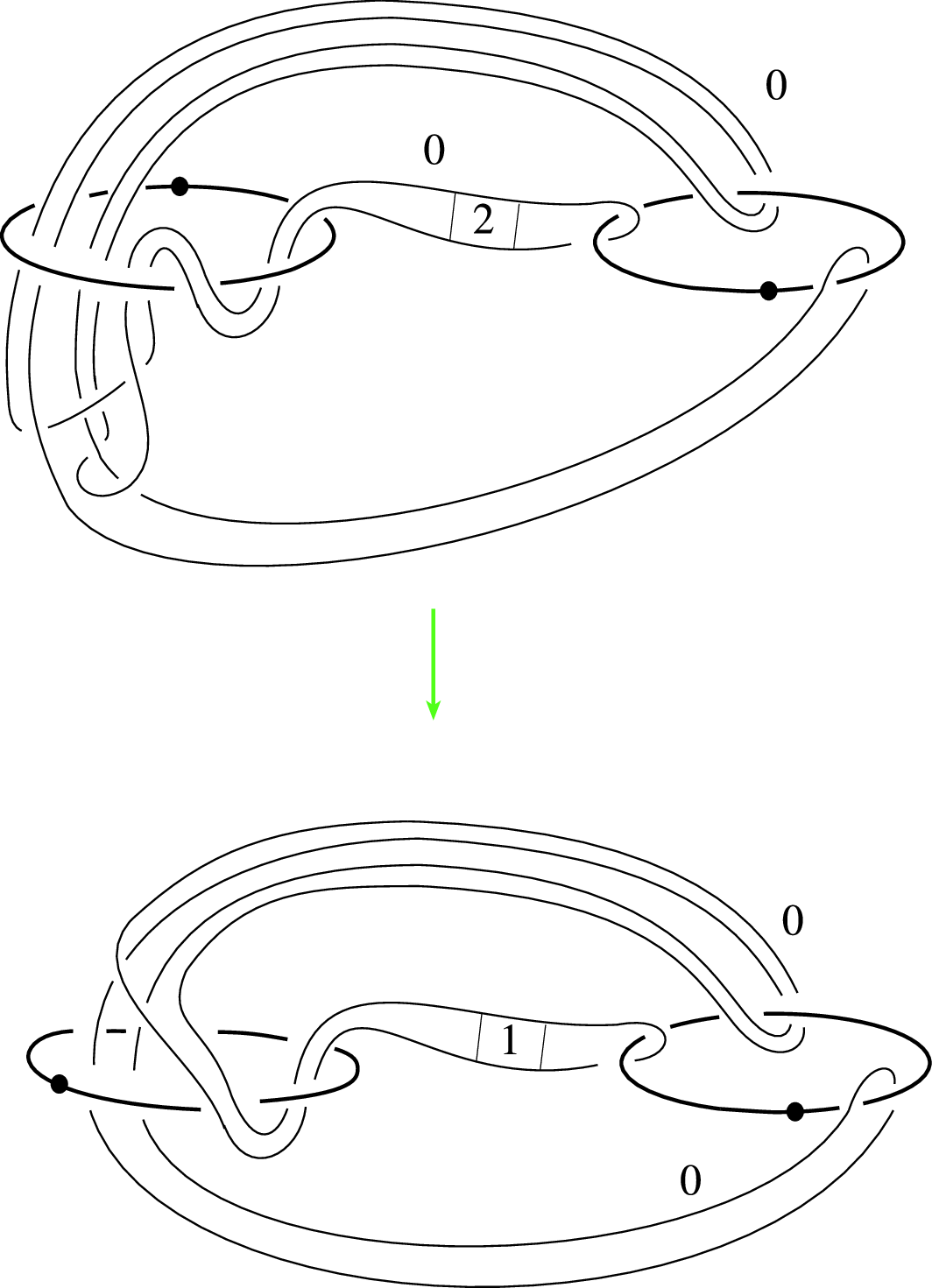}   
\caption{ } 
\end{center}
\end{figure}

    \begin{figure}[ht]  \begin{center}  
\includegraphics[width=.6\textwidth]{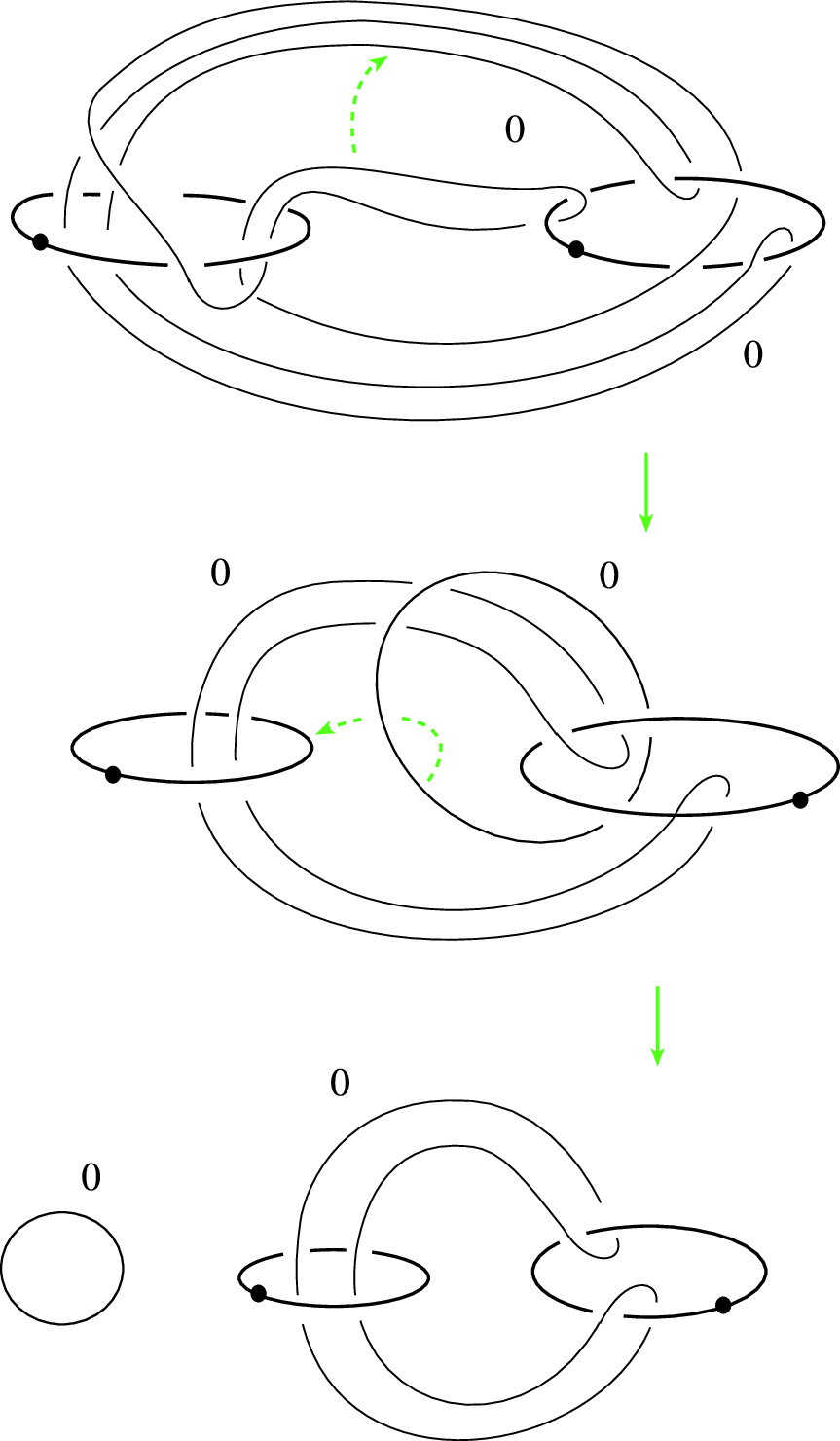}   
\caption{ } 
\end{center}
\end{figure}

    \begin{figure}[ht]  \begin{center}  
\includegraphics[width=.7\textwidth]{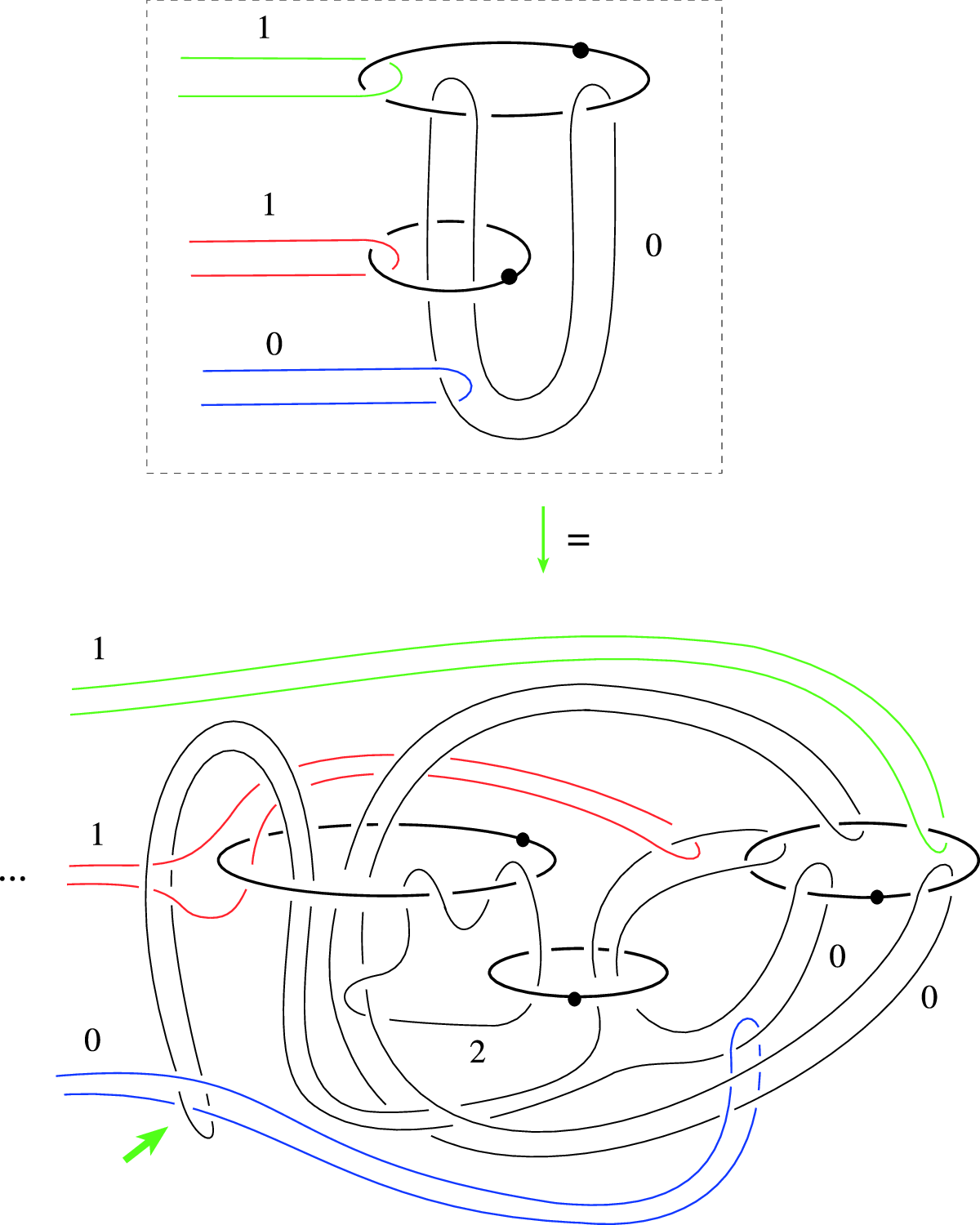}   
\caption{ $T^2 \times D^2$ containing a degree $2$ imbedded copy of itself, with some reference arcs drawn to  describe the diffeomorphism from the standard copy  $T^2\times D^2 \subset   S^{3}_{0}(K)\times S^1$} 
\end{center}
\end{figure}

    \begin{figure}[ht]  \begin{center}  
\includegraphics[width=.9\textwidth]{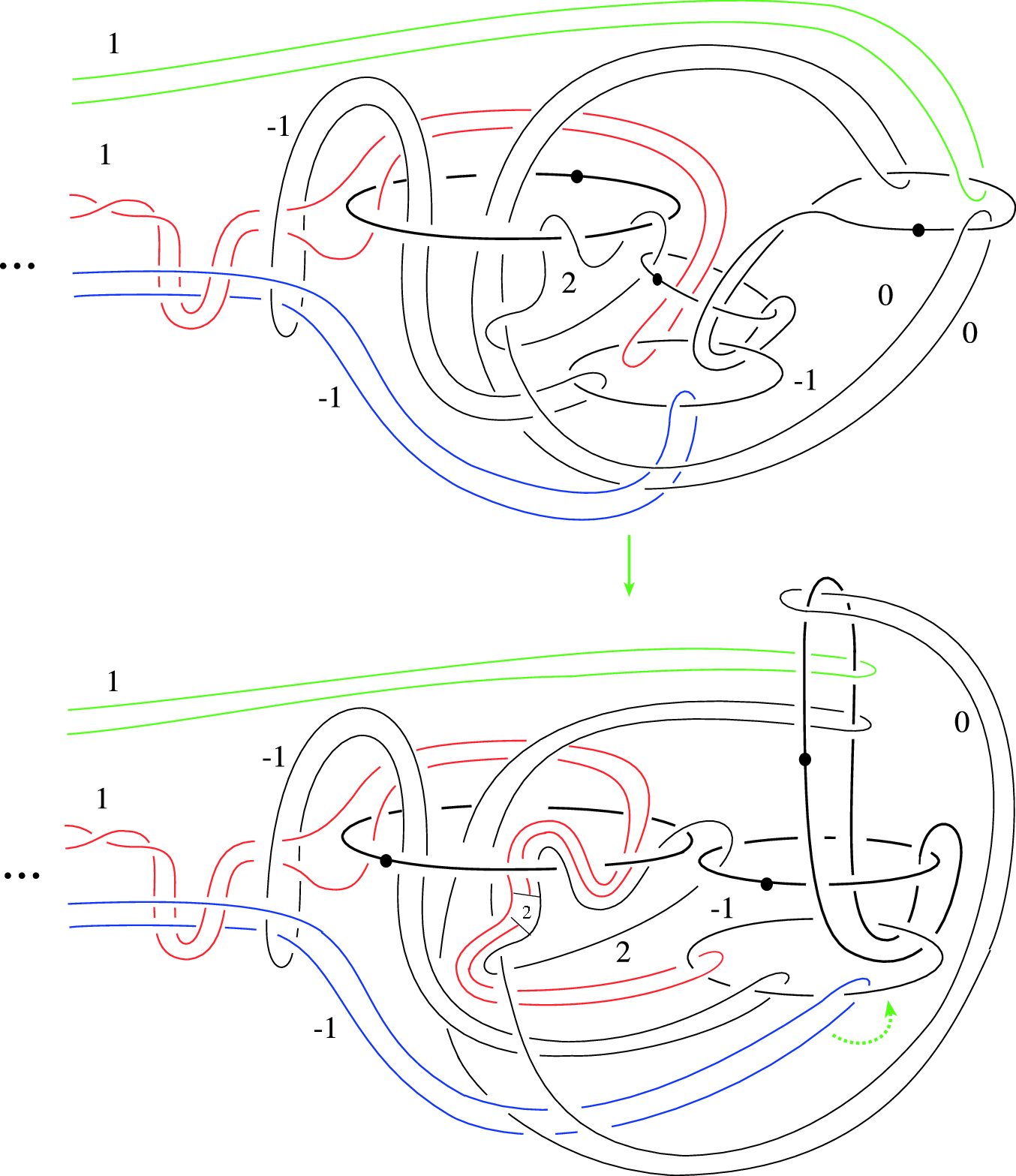}   
\caption{$(T^2\times D^2) \# \bar{CP}^2$ } 
\end{center}
\end{figure} 

    \begin{figure}[ht]  \begin{center}  
\includegraphics[width=.9\textwidth]{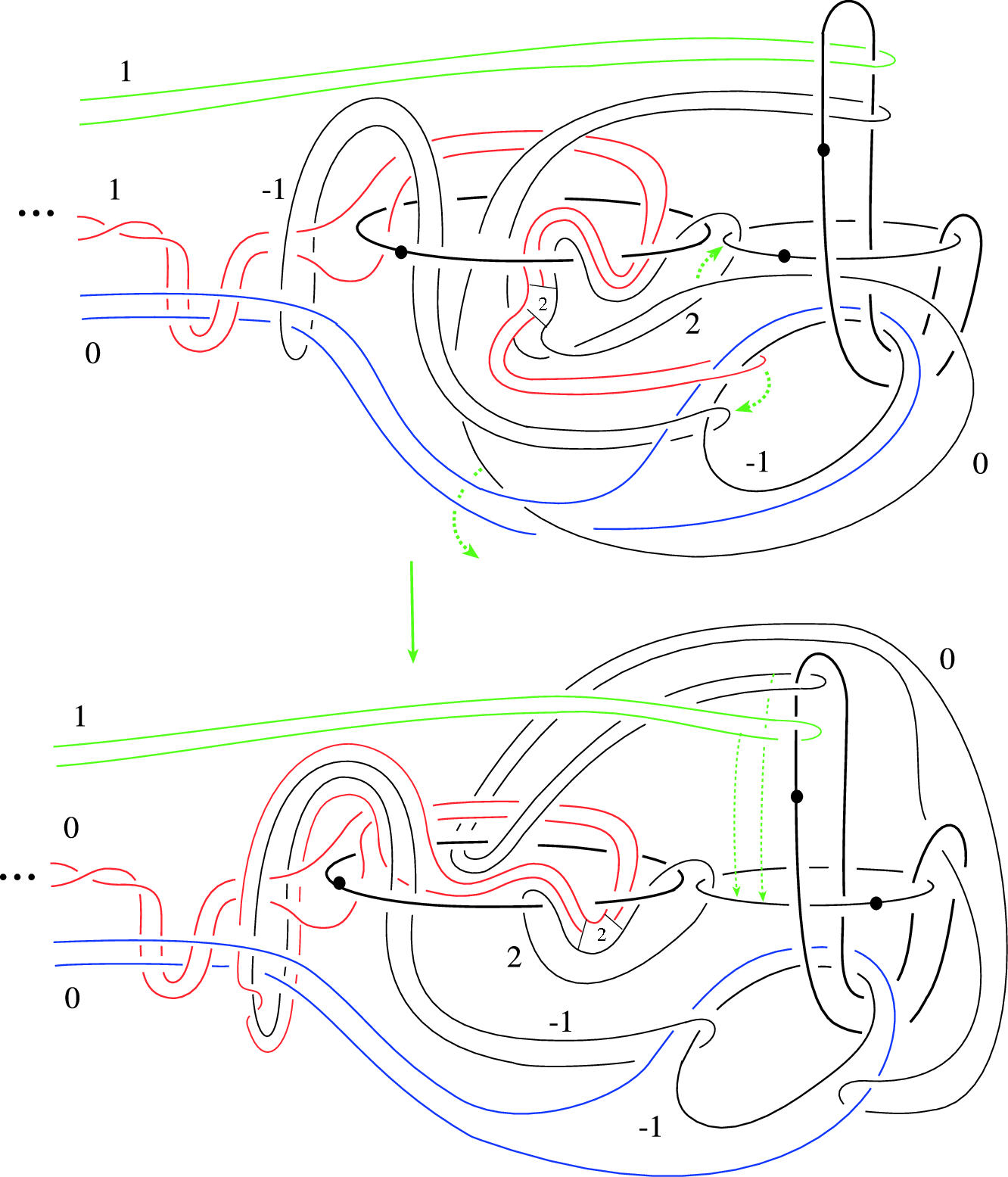}   
\caption{ $(T^2\times D^2) \# \bar{CP}^2$} 
\end{center}
\end{figure}

    \begin{figure}[ht]  \begin{center}  
\includegraphics[width=.9\textwidth]{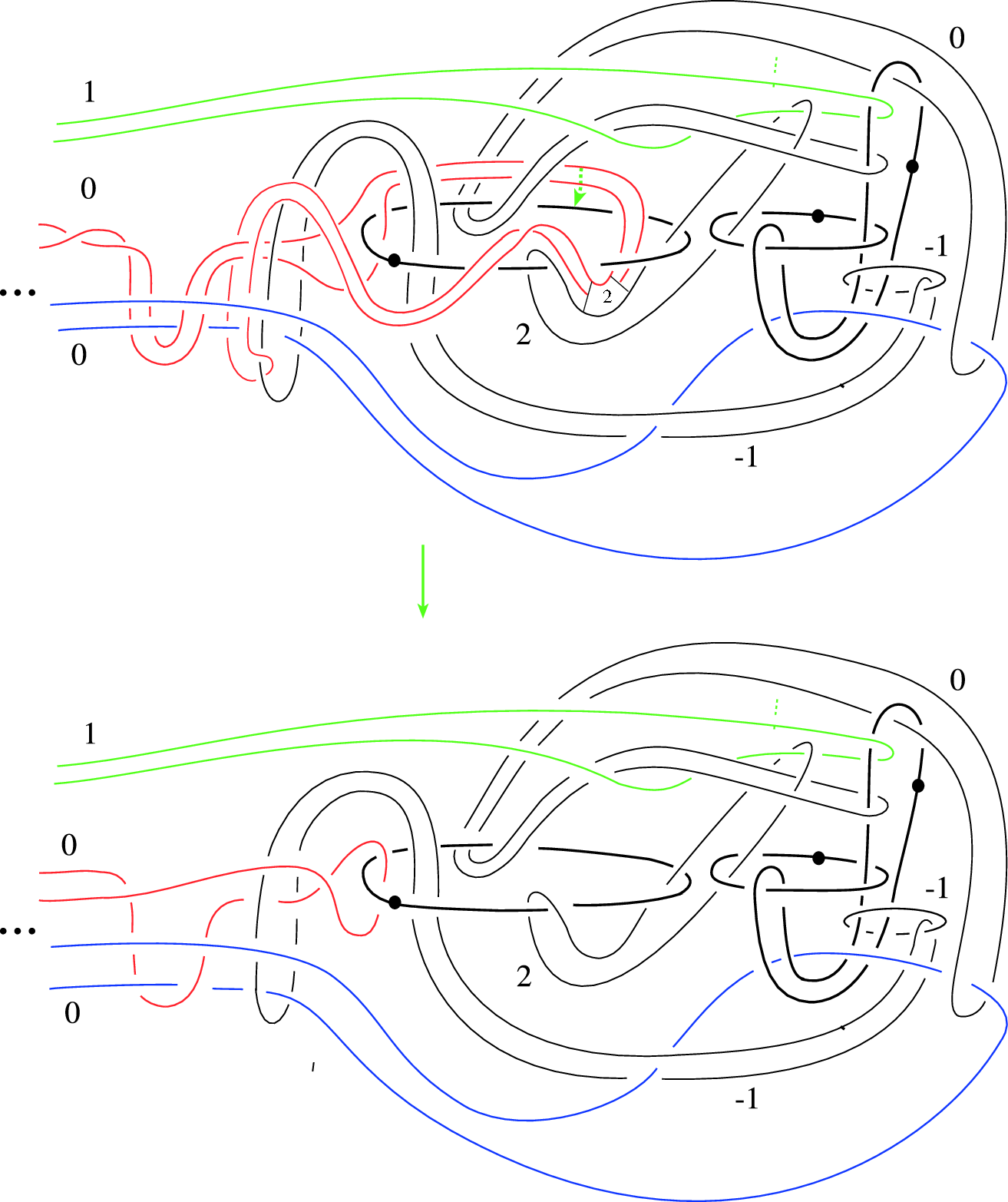}   
\caption{$(T^2\times D^2) \# \bar{CP}^2$ } 
\end{center}
\end{figure} 

    \begin{figure}[ht]  \begin{center}  
\includegraphics[width=.8\textwidth]{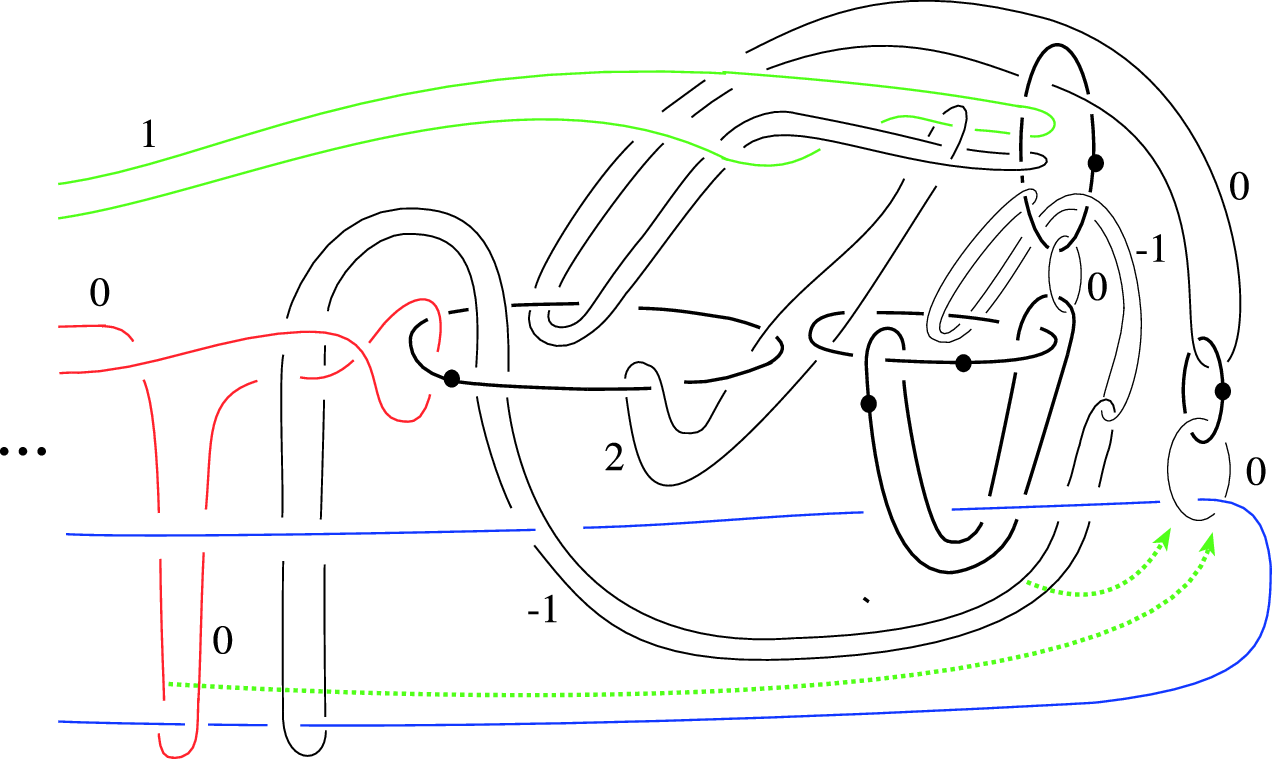}   
\caption{$(T^2\times D^2) \# \bar{CP}^2$ } 
\end{center}
\end{figure} 

    \begin{figure}[ht]  \begin{center}  
\includegraphics[width=.7\textwidth]{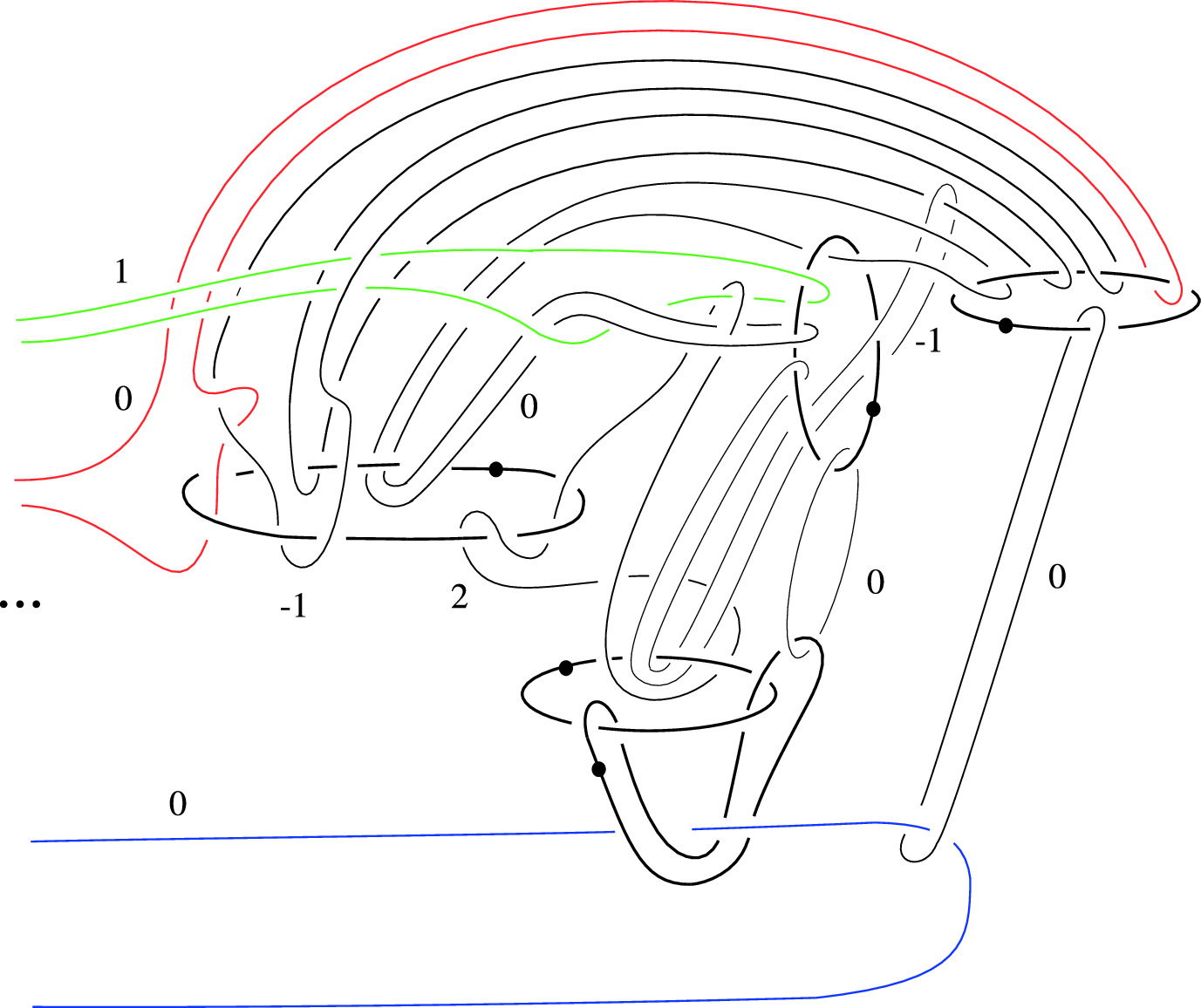}   
\caption{ $(T^2\times D^2) \# \bar{CP}^2$} 
\end{center}
\end{figure} 

    \begin{figure}[ht]  \begin{center}  
\includegraphics[width=1\textwidth]{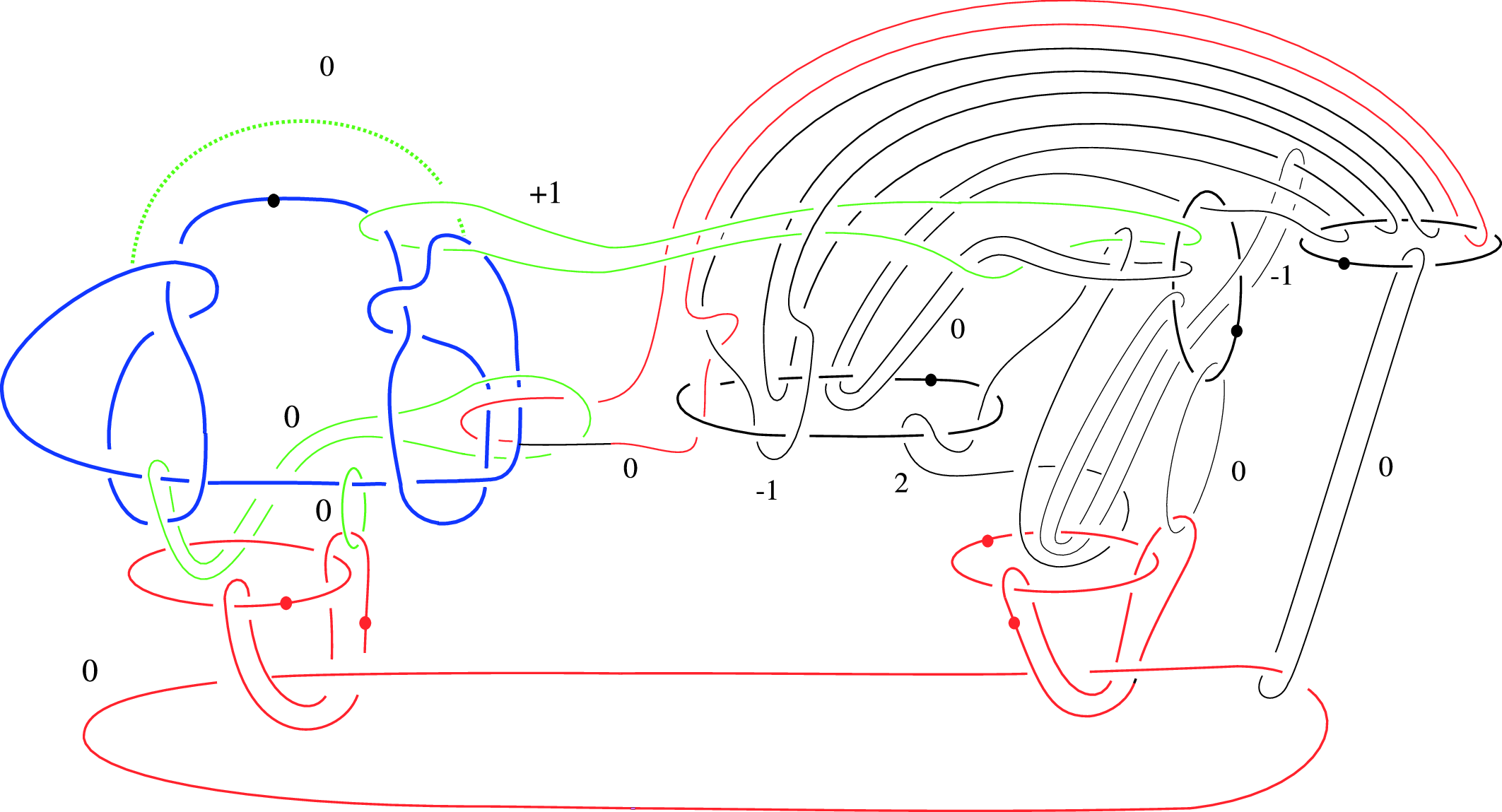}   
\caption{ } 
\end{center}
\end{figure} 

    \begin{figure}[ht]  \begin{center}  
\includegraphics[width=1\textwidth]{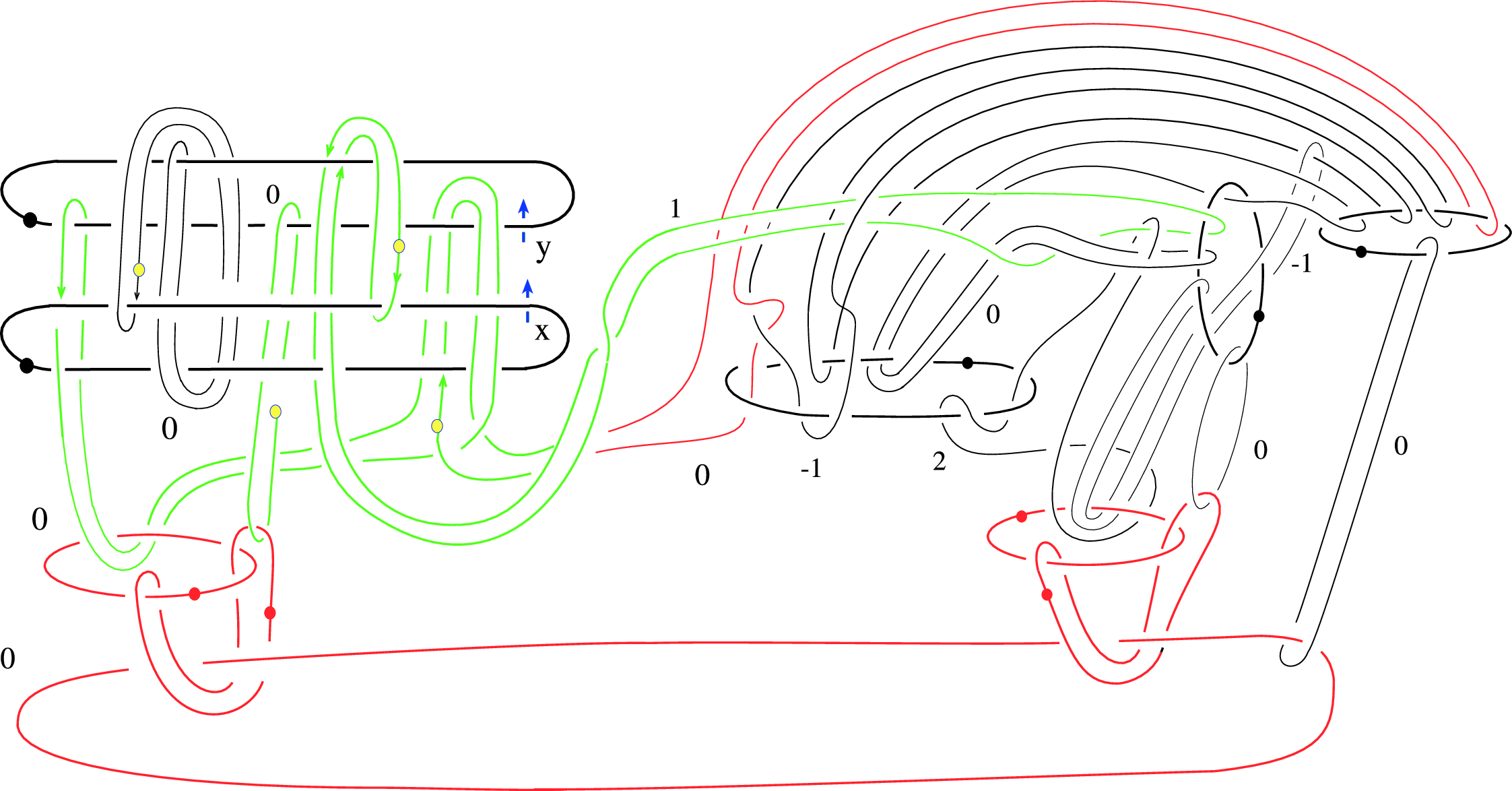}   
\caption{ } 
\end{center}
\end{figure} 

    \begin{figure}[ht]  \begin{center}  
\includegraphics[width=1\textwidth]{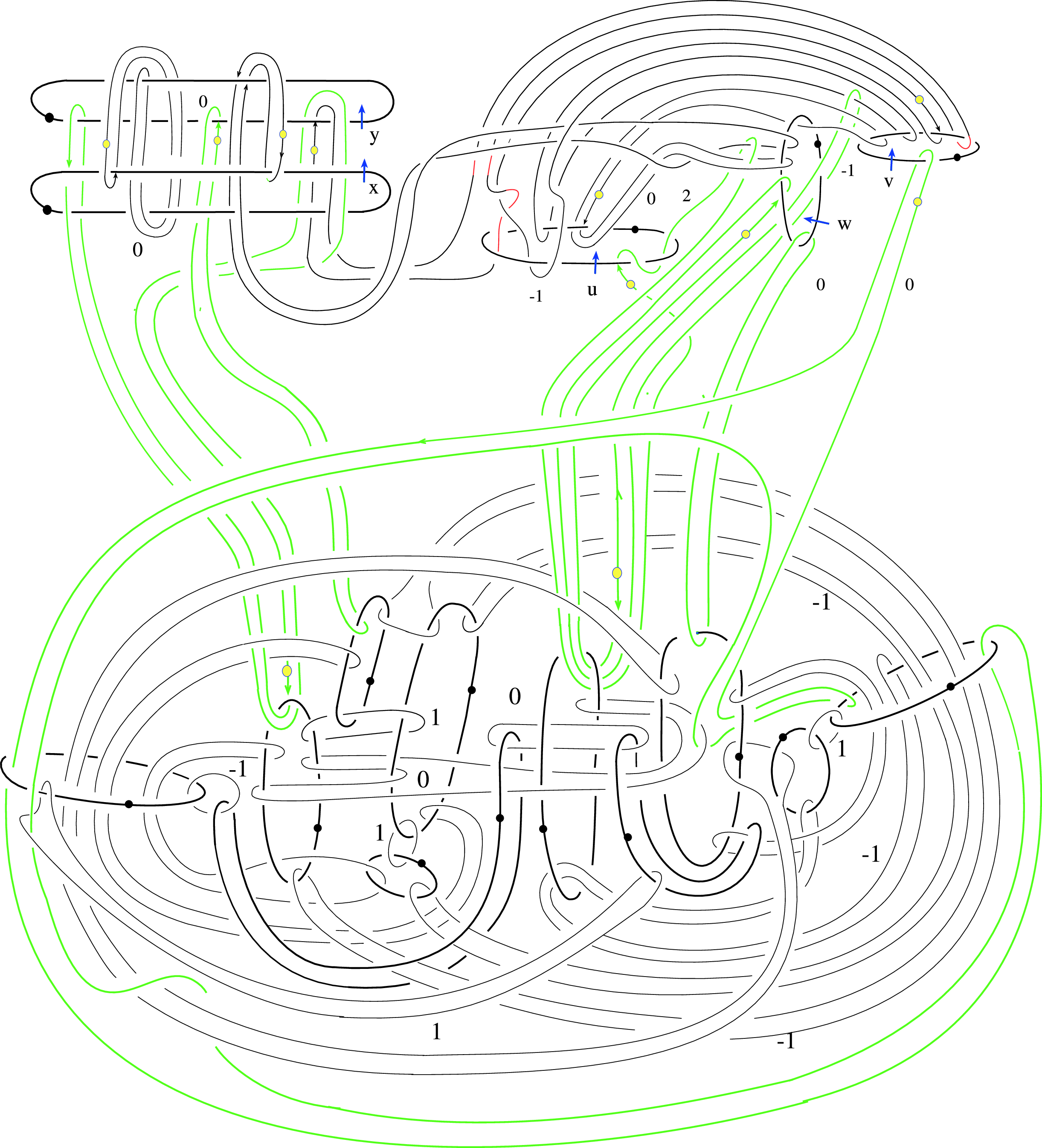}   
\caption{$N$ } 
\end{center}
\end{figure}

\end{document}